Cardano v Tartaglia:

The Great Feud Goes Supernatural

Tony Rothman**Words.**

Somewhere in Rome in October, 1570, the Brescian mathematician Niccolò

Tartaglia, "the stutterer," met with Aldo Cardano, son of Tartaglia's bitter enemy

Girolamo Cardano. In return for promises of gaining an appointment as a public

torturer and executioner, Aldo revealed to Tartaglia his father's whereabouts in

Bologna. Tartaglia hastened to that city and had Cardano arrested on charges of heresy for having cast a horoscope of Jesus Christ.

If you have heard this story, or some version of it, you are far from alone, for it is to be found in well-known books and on prominent websites. If you believe it, you are in good company as well, because the same forums pass the story on without qualification. If you believe it, though, you have been hoodwinked, because it is complete and utter nonsense.

* email: tonyrothman@gmail.com

The story of the epic sixteenth century feud between Girolamo Cardano and Niccolò Tartaglia over the solution to the cubic equation is justly one of the most famous in the history of mathematics. Its more colorful versions, involving the obligatory shifting alliances of the sixteenth century, subterfuge, betrayals and secret dossiers—let's not forget poison and syphilis—fairly scream for a theatrical presentation. Even a slightly sober investigation, however, reveals a less than Borgian scenario. The "cubic affair" in fact becomes a prime example of how scientific folktales, which have little or no basis in the historical record, nevertheless get passed up the great chain of existence until they become enthroned in the eighth heaven of print or cyberspace. In the case of the Great Feud, we are privileged to be able to trace the progress of the tale in an apparently straightforward manner. Nowhere in the strictly scholarly works on Cardano, for instance those by James Eckman,[1] Anthony Grafton[2] or Nancy Siraisi,[3] will one find any of the lurid details above, or indeed in the standard nineteenth century account, Henry Morley's loquacious two-tome biography of the

[1] James Eckman, *Jerome Cardan* (Johns Hopkins: Baltimore, 1946).

[2] Anthony Grafton, *Cardano's Cosmos* (Harvard: Cambridge, 1999).

[3] Nancy Siraisi, *The Clock and the Mirror* (Princeton: Princeton, 1997).

astrologer-physician.[4]   In the Italian literature on Tartaglia, the biographies of Masotti[5]

and Gabrieli[6] for example, any vengeful machinations—murderous or merely

injurious—are equally absent.  Even Oystein Ore's semi-popular work on Cardano,[7]

which though lacking references and unblushingly biased in Cardano's favor, more or

less adheres to known facts and avoids descent into sensationalism.   A discontinuity

occurs when one passes to the ultra-violet end of the spectrum.  There, aboard more

popular retellings, Hal Hellman's *Great Feuds in Mathematics*,[8] and Alan Wykes' *Doctor

Cardano, Physician Extraordinary*,[9] one decisively abandons the world of documents and evidence for realms unknown.

       That tabloid histories have supplanted mundane reality in numerous essays

suffixed by .edu  is perhaps less surprising than it is sad or amusing, depending on your momentary disposition.  Scientists, we must face facts, are suckers.  Beyond the hermetic world of scientific discourse, a significant percentage of folks fail to observe our clerical vows to facts, data, natural law and logic, and the same folks aren't above pulling a fast one.  It is not for nothing that magician James Randi in his investigations of pseudoscientific claims has steadfastly advocated that one needs someone trained in uncovering deception, not someone whose second nature presumes honesty.  Scientists, trusting souls, can be ruled out.  In their naivety they also, perhaps even more than most children, love a good story.  When a tale comes around that satisfies our analytical lust for the three C's: completeness, consistency and contingency ("no-plot-holes storytelling"), scientists' inherent gullibility leads us to accept it without question, in particular when it's too good to be true.

---

[4] Henry Morley, *Jerome Cardan* (Chapman and Hall: London, 1854).

[5] A. Masotti, *Niccolò Tartaglia*, in *Storia di Brescia*, **II**  pp.  587-617 (1963).

[6] Giovanni Battista Gabrieli, *Nicolò Tartaglia, Invenzioni, Disfide e Sfortune* (Brescia: 1986). This is the most complete account of Tartaglia's life I have found.

[7] Oystein Ore, *Cardano the Gambling Scholar* (Dover: New York, 1965).

[8] Hal Hellman, *Great Feuds in Mathematics: ten of the liveliest disputes ever* (John Wiley: Hoboken, 2006).

[9] Alan Wykes, *Doctor Cardano, Physician Extraordinary* (Frederick Muller: London, 1969).

**Absolute Truth (more or less).**

The tales surrounding the Great Feud do cry out loudly for a theatrical release; indeed my initial impetus to investigate them was to write a play about the episode, which—despite the stubborn intrusion of reality—I ultimately did, titling it *The Great Art*.  Much of the first half of what people believe they know about the famous affair is in fact true.[10]  The outstanding mathematical challenge of the early 1500s was to solve the cubic equation, in other words to find a "cubic formula" analogous to the famous quadratic formula, which had been known since antiquity.  By contrast, the cubic formula had eluded all attempts to find it and most mathematicians of the era, following Fra Luca, believed that such a solution was beyond the powers of human reason.

---

[10] See Ore's chap. 3 or the St. Andrew University MacTutor History of Mathematics website (henceforth MacTutor):
http://www-history.mcs.st-and.ac.uk/HistTopics/Tartaglia_v_Cardan.html

The Italian university system at the time curiously resembled our own, with tenure nonexistent and itinerant professors eking out an existence on temporary appointments. In such a milieu an important means of advancement were public "challenge matches," mathematical, medical and otherwise, which incidentally proved extremely popular with the citizenry. In 1535 a mathematician Antonio Maria Fiore challenged Niccolò Tartaglia (1499-1557) to such a contest. Mysteriously, Fiore had been boasting that he was in possession of the solution to the "depressed cubic," that is an equation of the form $x^3+ax = b$, where $a$ and $b$ are positive numbers. (At the time, the concept of a solution to the general cubic $ax^3+bx^2+cx+d = 0$, for any real coefficients, had yet to arise. Numbers reflected the positive physical world and hence negative numbers were highly suspect. The equation $ax^3+cx+d = 0$ was thus regarded as completely different from $ax^3+cx = d$, which in turn was completely different from $ax^3+bx^2 = d$. There were thirteen cases in all, which needed to be solved separately.)

Fiore posed to Tartaglia thirty problems, all of which boiled down to the depressed cubic. ("A man sells a sapphire for 500 ducats, making a profit of the cube root of his capital.  How much is the profit?"[11])  Several years earlier, Tartaglia had discovered how to solve the case $ax^3 + bx^2 = d$ and on the night of February 12-13, 1535 he perceived the solution to the depressed cubic as well.  Tartaglia was thus able to solve all of Fiore's problems within two hours and, for his own part, having posed problems that Fiore could not solve, easily won the match.  Tartaglia declined the thirty banquets that were the stakes of the contest.[12]

News of Tartaglia's victory spread throughout Italy and in 1539 Girolamo Cardano (1501-1576), who was preparing a book on mathematics, approached Tartaglia with a request for his solution.  After strenuous refusals Tartaglia finally relented when the two met in Cardano's house in Milan, on condition that Cardano never publish it. Cardano swore a sacred oath that he would not.  However, in 1543 he and his student Ludovico Ferrari (1522-1565) learned that Antonio Fiore had gotten the solution from his own teacher, Scipione del Ferro (1465-1526), who had discovered it three decades earlier, but never published.  Feeling released from his vow, Cardano published the solution, with considerably more praise for del Ferro than Tartaglia, as well as solutions to the other cases, in his 1545 book the *Ars Magna*,[13] which became the most important mathematical treatise of the sixteenth century.

---

[11] *The History of Mathematics: A Reader*, John Fauvel and Jeremy Gray, eds. (The Open University: London, 1987), p. 254.

[12] Gabrieli gives numerous excerpts from the Tartaglia-Cardano dispute.  For English-language excerpts see Fauvel and Grey, pp. 254-256 ; MacTutor and Ore.

At that point, one might say without exaggeration that all hell broke loose. Tartaglia, in his own book *Quesiti et Invenzioni Diverse* (*Various Questions and Inventions*) of 1546, accused Cardano of theft and a violation of a sacred trust (or perhaps a financial one—challenge matches were after all worth good money). Cardano, by then Italy's most famous physician and astrologer, evidently did not want to enter into a public dispute with Tartaglia and turned over the matter to Ferrari, who very publicly challenged Tartaglia to a contest. Vicious manifestos flew back and forth between the two for eighteen months. "You make up proofs in your own head and thus they usually have no conclusion." "I truly do not know of any greater infamy than to break an oath, and this holds not only in our own, but in any other religion." "With these lies you attempt to convince the ignorant that your statements are true." " I honestly expect to soak the heads of both of you in one fell swoop, something that no barber in all Italy can do." "You are a devil of a man, wanting to be an inventor when you have the head of an adder, which can understand nothing."[14]

---

[13] Girolamo Cardano, *Ars Magna, or The Rules of Algebra*, translated by T. Richard Witmer (Dover: New York, 1993).

[14] Facsimiles of the original manifestos, first published by Enrico Giordani in 1876, are now available online in Latin and old Italian at
http://books.google.com/books?id=rBB1mTyvRDsC&source=gbs

Apart from reputation, under dispute were thirty-one questions each combatant had proposed to the other on algebra, geometry and philosophy. On August 10, 1548, the two antagonists and a large crowd of Ferrari's supporters met at The Church in the Garden of the Frati Zoccolanti in Milan for the final judging, presided over by the governor of Milan himself. No record exists of what exactly transpired during the occasion. It is generally accepted that Ferrari won, because Tartaglia slipped away during the first night, although from Niccolò's recollections one might conclude that he couldn't get a word in edgewise.

**Disputable.**

All of this is fairly well documented: the *cartelli* and problems exchanged between Tartaglia and Ferrari exist, and in his books Tartaglia gives verbatim accounts of his letters and meeting with Cardano.[15] There is little reason to suspect that Tartaglia's version is far from the truth: Apart from the fact that Niccolò appears to have been a pack-rat until the end of his life, Cardano in effect never disputed his claim

---

For excerpts in Italian, see Gabrieli, or Luigi di Pasquale, "I cartelli di matematica disfida di Ludovico Ferrari e i controcartelli di Nicolò Tartaglia," I, *Period., Mat.* (4) **35,** 253-278 (1957); II, *Period., Mat.* (4) **36** 175-198 (1957). For English excerpts see Ore.

[15] See sources already cited.

in the *Ars Magna*.   Some authors, for example Witmer[16] and Hellman,[17] argue that Ferrari (who was present at the meeting between Cardano and Tartaglia) later vociferously denied that Cardano had ever sworn such an oath.  I find no evidence that this is the case.  The relevant passage is from Ferrari's second *cartello*:

> First of all let me remind you, so that you don't remain astonished and wonder where I have heard all your lies, as if by a revelation of Apollo, that I was present in the house when Cardano offered you hospitality and I attended your conversations, which delighted me greatly.  It was then that Cardano obtained from you this bit of a discovery of yours about the cube and the *cosa* equal to a number,[*] and this languishing little plant he recalled to life from near death by transplanting it in his book, explaining it clearly and learnedly, producing for it the greatest, the most fertile and most suitable place for growth. And he proclaimed you the inventor and recalled that it was you who communicated it when requested.

---

[16] See Witmer's preface to the *Ars Magna* [13], p. xviii, note 26.

[17] Hellman, p. 18.

[*] In sixteenth century Italy, the unknown was referred to as the *cosa* (the thing).  "The *cosa* and the cube equal to a number" was therefore the expression for the depressed cubic $x^3+ax = b$.

What more do you want? "I don't want it divulged," you say. And why? "So that no one else shall profit from my invention. "And therein, although it is a matter of small importance, almost of no utility, you show yourself un-Christian and malicious, almost worthy of being banned from human society. Really, since we are born not for ourselves only but for the benefit of our native land and the whole human race, and when you possess within yourself something good, why don't you want to let others share it? You say: "I intended to publish it, but in my own book." And who forbids it? Perhaps it is because you have not solved it entirely….[18]

Polemics one sees in abundance; an oath or its denial, no. Robert Kaster of Princeton University has graciously checked the facsimile of the entire Latin original for me and finds no mention of the oath elsewhere. Nor does Ore, who presents this translation, claim any denial of oath on Ferrari's part. To all appearances it is merely Ferrari's justification, on the part of the human race, for Cardano's publication of the cubic formula, and his admonition to Tartaglia to stop kvetching.

Hellman also gives credence to Alan Wykes' claim in his book *Doctor Cardano, Physician Extraordinary* that Cardano in fact worked out the formula for himself and then by "a slip of pen or memory, he wrote that Tartaglia had communicated the discovery to him and given him permission to use it."[19]  Wykes, in a manner that will become familiar, gives no justification for this fabulous assertion, which requires that Tartaglia invented not only the meeting between himself and Cardano, but their entire correspondence.  It also makes Ferrari's eyewitness account impossible.  For the record, in the *Ars Magna* Cardano writes, "[Tartaglia] gave [the rule] to me in response to my entreaties, though withholding the demonstration."[20]

---

[18] With some minor corrections from Robert Kaster this is the translation given by Ore, p. 94.

Oath aside, to this day the larger discussion centers on whether Cardano's actions were justified, given that Tartaglia had failed to publish his results in the decade after his contest with Fiore.  I intend to avoid that particular debate.   For divers opinions the reader may want to see Eckman's detailed study.[21]  (On the matter of the oath, Eckman writes, "There is, of course, no doubt as to the breach of faith on the part of Cardan.  It was flagrant, even if allowance is made for the moralities of the sixteenth century in respect to mutual relationships."[22])  I do point out that statements, beginning with Ferrari's, to the effect that Tartaglia stood against the progress of science by intending to keep his discovery secret, appear grounded less in reality than in rhetoric.

---

[19] Wykes, p. 115.

[20] Cardano, *Ars Magna*, p. 96.

[21] Eckman, chap. 4.

[22] Ibid., p. 64.

Even in the midst of his diatribe Ludovico recollects that Niccolò had protested only that he wanted to publish it himself.   To be sure, in 1539 Tartaglia had said to the bookseller Zuan Antonio de Bassano, who acted as intermediary between himself and Cardano,  "Tell *Eccellenza* that he must pardon me: when I propose to publish my invention, I will publish it in a work of my own, not in the work of another man, so that *Eccellenza* must hold me excused."[23]

Of course Tartaglia did not publish; nevertheless his excuse was evidently plausible: for many years he was occupied with the first translation of Euclid into any living language (Italian, 1543), and a modern edition of Archimedes (1544), both signal events in the history of mathematics.   Indeed, in 1541 he wrote to his English pupil Richard Wentworth, assuring him that he would publish his formula once these works were done.[24]   Tartaglia may have also lost his entire family at about the same time.[25] And then Cardano beat him to the punch.

As it turns out, a year after the publication of the *Ars Magna*, Tartaglia published his *Quesiti*, where one finds this striking passage in the dedication:

---

[23] Ore, p. 66, and MacTutor.

[24] *Mechanics in Sixteenth Century Italy*, translated and annotated by Stillman Drake and I.E. Drabkin (University of Wisconsin: Madison, 1969);
*Metallurgy, Ballistics and Epistemic Instruments,* The Nova Scientia *of Nicolò Tartaglia, a new edition*, Matteo Valleriani et al. eds. (Edition Open Access: Berlin: 2013):
 http://www.edition-open-access.de/sources/6/index.html

[25] Drake and Drabkin, p. 21; Gabrieli, p. 20.

I reflected that no small blame is attached to that man who, either through

science, his own industry or through luck, discovers some noteworthy thing but

wants to be its sole possessor; for, if all our ancients had done the same, we

should be little different from the irrational animals now. In order not to incur

that censure, I have decided to publish these questions and inventions of mine.[26]

Unless one believes that this statement was forced by publication of the *Ars Magna*, it does not appear to be of a man unwilling to divulge his results.

Regarding the general view, implicit in Ore's work, that Tartaglia's position made him the last "medieval man" who put personal gain over communal progress, one might at this juncture bemoan the fact that those writing about the feud have been mathematicians rather than physicists. Tartaglia's first book, the *Nova Scientia* of 1537,

was in fact the earliest attempt to treat the trajectory of projectiles by mathematical

means, and it surely provided the model for Galileo's later *Two New Sciences.* In the

*Quesiti*, Tartaglia became probably the first natural philosopher to openly challenge

Aristotelian mechanics. It is interesting that while Cardano's publication of the cubic formula resonates with today's "open source" culture, Tartaglia's reasons for hesitating to publish his results on ballistics ("it was a blameworthy thing…a damnable exercise, destroyer of the human species…[and] I burned all my calculations…[27]) might have been written by today's anti-nuclear movement.  Only under threat of a Turkish invasion did Tartaglia change his mind.  It is also curious that the St. Andrew University MacTutor History of Mathematics website, which is fairly comprehensive, does not even mention Tartaglia's major work, the *Trattato Generale di Numeri et Misure* of 1556, usually considered one of the most important textbooks on arithmetic of the sixteenth century.

---

[26] Drake and Drabkin, p. 99.

**Falsifiable.**

If scholarly DNA requires arguments about everything, one thing is fairly impervious to even academic genetic coding:  After the face-off between Tartaglia and Ferrari in 1548, the historical record rapidly grows mute.  Little is factually known about Tartaglia's life apart from the occasional public document and autobiographical passages scattered throughout his mathematical works.  As we know, however, Nature abhors a vacuum, and it may well be the vacuum that has inspired authors to fill it with tales that extend the feud to literally the supernatural domain.

It is true, as Ore relates, that after the misadventure with Ferrari, the patrons who in early 1548 had invited Tartaglia to his native Brescia to lecture on Euclid did an about-face and refused to pay him for his labors. Niccolò lost eighteen months' salary and was forced to return to Venice, where had had lived since 1534, and continue his livelihood as a private mathematics teacher. But as plausible as it might seem that his hosts' bad faith was the result of his poor showing in Milan[28]—contingency, after all— there is no documentary evidence that this is the case. In fact, Tartaglia continued to lecture in Brescia for another year after the historic showdown. For this reason Gabrieli argues that the two events are unconnected.[29]

---

[27] Ibid., p. 68.

There can't be any doubt that Tartaglia remained extremely bitter about what had transpired and even in his last work, the *Trattato Generale*, he returned to the problems posed a decade earlier in the manifestos, making scornful remarks about his opponents' solutions. Nevertheless, all stories—*all*—that Tartaglia devoted the remainder of his life to revenging himself against his nemesis are apocryphal, in the original sense of the word, or plainly false. The most recent retelling is Hal Hellman's 2006 *Great Feuds in Mathematics*,[30] already mentioned, which I now quote at length because it provides a concise compendium of what have become the standard rumors and legends surrounding the Cardano-Tartaglia affair. By the mid sixteenth century,

---

[28] Ore, p. 105.

[29] Gabrieli, p. 85.

[30] Hellman, pp. 23-24.

the Roman Inquisition and Counter Reformation were underway.   In the decades after

the Ferrari-Tartaglia contest

…Scholars of all sorts were under suspicion, but somehow Tartaglia had

managed to place himself satisfactorily.  Cardano could find no employment

and, according to Wykes, "it was Tartaglia who was the instigator of most of the

refusals that met him in College and University.  It was simple enough, with the

network of the Inquisition flourishing in city, vineyard, village and public

square, to keep a shadowy hand on the shoulder of any citizen, great of small."

This was just the warmup, though.  On October 13, 1570, almost a quarter of a century after publication of *Ars Magna*, Tartaglia served up a double blow.

Using Cardano's own son Aldo as in informant as to Cardano's whereabouts,

Tartaglia handed him to the Inquisition.  Tartaglia had been collecting evidence against Cardano for years.  Among this "evidence" was Cardano's rejection of the pope's invitation that he become the pope's astrologer and physician.

Tartaglia pointed to the "sarcasm" evident in Cardano's comment that "His

Holiness by his study of astrology has surely raised himself among the greatest

of such scientists and has no need of help from such as myself."

Cardano's horoscope of the life of Jesus was also damning, as were a variety of other statements that, taken out of contexts, could be construed as blasphemous.  In one of his publications, for example, he had suggested that God is a universal spirit whose benevolence is not restricted to holders of the Christian faith.  Today he might be admired for such an ecumenical statement; at the time it was apparently a dangerous idea.

And so it went.  Cardano, fortunately, was not subjected to torture or put to death, but he was thrown into jail.  He sought desperately for help and was able to reach out to an official in the church, Archbishop Hamilton, who had in the past asked to be called upon if need be.  The archbishop came through for Cardano, who was released a few months later.  It was just in time, for not long after, the archbishop's own fortunes changed; he was captured by the forces of Mary, Queen of Scots, and beheaded.

Tartaglia finally had had his revenge.  Cardano lived on in obscurity in Rome, where he worked on his autobiography, which is one of the works that has come down to us in full.  He probably never knew, and just as well, that his daughter Chiara had died of syphilis, and that it was Aldo who betrayed him to the Inquisition and who was rewarded with an appointment as official torturer and executioner in Bologna.

Cardano died on September 20, 1576.  Less than a year later, Tartaglia followed him to the grave.

As signaled in the introduction, the same stories, that "Cardano himself was accused of heresy in 1570 because he had computed and published the horoscope of Jesus in 1554," and that "apparently, his own son [Aldo] contributed to the prosecution, bribed by Tartaglia," can be found in the Wikipedia entry on Cardano.[31]  The contention that his daughter Chiara died of syphilis is so widely spread on the Internet that specific references are unnecessary.  According to one essay, the tragedy prompted Cardano to write one of the earliest treatises on the disease.

What truth to these tales?  First, the contention that Cardano was unable to find a job, while Tartaglia "managed to place himself satisfactorily," is completely counterfactual.  The 1550s saw Cardano at the height of his fame, with a professorship in Pavia, at least one genuine bestseller (*De Subtilitate*) and invitations by European potentates (e.g. Archbishop Hamilton of Scotland, whom Cardano cured of asthma to great acclaim).  Throughout the 1560s, Cardano remained relatively prosperous, although he resigned from the University of Pavia, evidently because of accusations of pedophilia,[32] and moved to a lectureship at the University of Bologna.  Tartaglia, on the other hand, returned to Venice in poverty and remained desperately poor until his death, bequeathing only books to his publisher, brother and sister, as well as a few household items to the latter. [33] (Niccolò had hardly been interred before the publisher, alas, made off with all the books.)  Cardano was indeed arrested in Bologna on October 13, 1570 for impiety,* although nowhere in his writings does he disclose the reasons, and no records of the proceedings have come to light.  It is possible that his arrest was due to his 1554 horoscope of Jesus Christ, but this has never been established, even if it is consistent.[34]  He was released from prison after three months due to intervention by his friends Cardinals Morone and Borromeo (not Hamilton), held under house arrest for a time, then invited to Rome, where he spent the last five years of his life provided for by the pope, continuing to practice medicine but no longer allowed to teach or publish.

---

[31] http://en.wikipedia.org/wiki/Gerolamo_Cardano. (Note added: Since posting of the first version of this preprint, the erroneous information Wikipedia entry has been removed.)

[32] Grafton, p. 188.  The matter is also discussed obliquely in Cardano's *Book of My Life*, pp. 96-99; see note 40, below.

As for the remaining stories, tales of Aldo Cardano's complicity in his father's arrest have been alternately surfacing and submerging since at least the nineteenth century,[35] but Tartaglia's role in the affair has been most notably propagated, if not altogether invented, by Alan Wykes, whose *Doctor Cardano, Physician Extraordinary* Hellman follows.  In Wykes' account not only can Tartaglia, a poor Venetian mathematician ("in whom the seeds of instability had been nourished by childhood

---

[33] Gabrieli, pp. 104-110.

* It is often said that Cardano was arrested on the charge of "heresy," but again, the exact charges have never come to light and it should be borne in mind that the Church was careful to distinguish among various religious crimes.  For example, Galileo was charged with "suspicion of heresy," a lesser offense than heresy.

[34] See Eckman, p. 33 *et seq.*.

[35] Ibid., pp. 32-33.

environment and had grown into weeds choking the flowers of his own brilliance")

worm his way into the good graces of the governor of Milan in order to thwart

Cardano's advancement, but in doing so he is able subvert his enemy by disclosing the horoscope of Jesus to a papal emissary.[36]  Wykes' plotting is impressive and lurid.   Too impressive, too lurid.  I have reluctantly come to the conclusion that his work was either written from memory without double-checking sources, or is a deliberate literary hoax.

One should of course think twice before imputing motive, but were Wykes alive, I

would certainly ask him to explain himself.

It is easiest to deal with Wykes' book by beginning at the end.  The closing is:

"[Cardano] died on 20th September 1576, a man not without greatness in an age of great

and cruel men.  Less than a year later his enemy Tartaglia died also."[37]

In fact, Tartaglia died on the night of 13-14 December, 1557, nineteen years before

Cardano.  This is not a matter of conjecture or debate: his Last Will and Testament exists

and has been published; Tartaglia was buried in the church of San Silvestro in Venice according to his wishes.[38]  Of course, Wykes' error makes most of the above claims impossible, by chronological protection.  It is nevertheless instructive to see how he justifies, for example, the tale that in 1570 Tartaglia bribed Cardano's son Aldo into turning his father over to the Inquisition.  Wykes' writes:

> The boy Aldo, to whom I had promised the reward of the appointment of public torturer and executioner in that city [Bologna], came to me in Rome with the intelligence that his father was in Bologna, awaiting an interview with the syndics.  I thought to myself, 'Ah!  This will be pleasant, to raise his hopes that at last the restrictions are about to be lifted from him and then, an instant before the realization of those hopes to cast him into prison.  And so it was.  I hastened to Bologna, and there he is still sheltered, in the ruins of a hovel, awaiting an ascent to his former status.  I instructed the guards to arrest him as he set out for his appointment.[39]

---

[36] Wykes, p. 117, pp. 120-121.

[37] Ibid., p. 176.

[38] Gabrieli, pp. 104-110.

The context makes clear that Wykes intends the reader to believe this passage was written by Tartaglia.  The fact that in 1570 Tartaglia had been dead for thirteen years should be sufficient reason to doubt it.  Additionally, there is no evidence that Niccolò was ever in Rome, Cardano's presence in Bologna was hardly a secret, and his presence in a hovel, really?—he had been awarded the high honor of being made a citizen of the city.  Who then wrote the passage?  In Wykes' book, it is tagged "footnote

---

[39] Wykes, p. 174.

2" for chapter eighteen, but in the endnotes for that chapter, a source for footnote 1 is listed and nothing more.   Reference 2 is simply missing.   Given that I have found it virtually impossible to confirm a single citation in Wykes' book, I would not be surprised if he invented it himself.[*]

Here I must turn to Cardano's autobiography, *De Vita Propria Liber*, or *The Book of My Life*,[40] which is one of the Renaissance's most famous memoirs and the work through which we know most about the author.   In it Cardano is remarkably frank about his failures as a father and the disasters of his two sons, the elder Giambattista (1534-1560) who was executed for poisoning an adulterous wife, the younger Aldo (1543-?), who was arrested on numerous occasions for theft.   After Aldo burglarized his father's own home in 1569, Cardano had him imprisoned and disinherited him.   Wykes makes extended assertions[41] that Aldo acted as a torturer and executioner and that Cardano knew it via public accounts ("Messer Aldo Cardano, executioner, for torturing by rack and vice, Valentino Zuccaro, 3 scudi."), but nowhere in *The Book of My Life* does

---

[*] The inadequate citations throughout Wykes' book make it extremely difficult to verify anything.  The few citations that are given are to titles only and never include page numbers.  Ore also fails to give references and certain quotations appear to me dubious (e.g., the unending adjectival string by which Cardano describes his own character on Ore's p. 25 is not to be found in Cardano's *Book of My Life* (next footnote)). When I have been able to track down others, however, they appear reasonably accurate.

[40] Girolamo Cardano, *The Book of My Life*, translated by Jean Stoner (New York Review Books: New York, 2002).  This translation originally appeared in 1929.

[41] Wykes, pp. 151-152.

Cardano mention any such activities.  The only source Wykes offers for his claims is

Cardano's *De Consolatione*, which was published when Aldo was negative one year old.[42]  A precocious child indeed.

A similar haze surrounds Cardano's daughter, Chiara (1536-?).  Wykes writes

that by the age of sixteen Chiara had seduced her elder brother Giambattista.[43]  No

reference is given.  He does present a single-sentence quotation "There was nought of

honesty at all in her whoring," which points us to Peter Martyr Vermigli's *Loci*

*Communes*.[44]  He next gives an extended excerpt from a letter by Chiara's husband

Bartolomeo Sacco, in which Sacco writes, "Not only have you shed upon me the great

pox in the person of your unclean daughter, but you have given me a wife whose

demands night and day are more than can be met by the staunchest lover of couch

pleasures…"[45]  The missive becomes far more graphic, ending with the husband's threat to seek an annulment of the marriage.   Again, no reference is given.   Wykes does provide a source for two subsequent passages regarding Chiara: "A young woman still, she was brought to book of the Spanish disease and her own sad flux."  Chiara's sterility was due to her incestuous relationship with Giambattista and "the exaction of the price" by the ecclesiastical courts for this crime "was endless."  The citation for the quoted passages is Cardano's *Book of My Life*.

---

[42] http://books.google.com/books/about/De_Consolatione.html?id=evs5AAAAcAAJ

[43] Wykes, p. 142.

[44] Ibid.

[45] Ibid., p. 149.

What are we to make of all this?  One might scratch one's head for a moment to ask why Peter Vermigli, a famous Florentine theologian, would be writing about Chiara Cardano, yet alone in a compendium of theological practices.  In answering this question I am limited by my inability to read Latin, but I have checked all the English translations of Vermigli's works at Princeton University (which do not include the *Loci Communes*) and there are only two passing references to Girolamo Cardano and none to Chiara.  The *Loci Communes* itself is now available online as a Google Book.[46]  Its index contains no mention of Cardano or his daughter.[*]

As for Wykes' references to *The Book of My Life*, we are immediately confronted by Cardano's own statement, "From my daughter alone have I suffered no vexations beyond the getting together of her dowry, but this obligation to her I discharged, as was right, with pleasure."[47]  In fact, I challenge anyone to find the passages Wykes cites in *The Book of My Life*.  Initially, I assumed that he must have worked from a more complete edition, but in his bibliography Wykes lists the translation he used as the one

---

by Jean Stover [sic].  The 1929 translation by Jean Stoner is the only one into English that I am aware of.  Under normal circumstances I would assume this was a simple misprint; in light of the rest…Neither does Morley in his biography of Cardano mention any such behavior on Chiara's part.

In a word, I have found only one "documented" contention apart from Wykes' that Chiara Cardano ended her life as a prostitute or died of syphilis.  Eckman does cite H. Kümmel as writing in 1910, "*Eine Tochter, das einzige Kind, das ihm geblieben war, brannte mit einem Galan durch und endete als Dirne*,"[48] or, "A daughter, the only child left to him, eloped with a gallant and ended up being a prostitute."  On the other hand, Eckman himself says of this passage that he knows of no authority for it.  Given that Chiara married Bartolomeo Sacco, a patrician, almost certainly while Girolamo's sons were still with him, it is difficult to see what authority there could be.

To sum up, as far as I am able to determine, all the direct quotations in Wykes' book from family and household members are either loose paraphrases from Cardano's *Book of My Life* or fabrications.  And incidentally, Cardano neither invented, nor claimed to invent the universal joint, or Cardan shaft—another popular pass-me-down that can be found in Wykes' book[49] and on Wikipedia—but only a chair that could be kept level on an incline.[50]

---

[48] Eckman, p. 33.

[49] Wykes, p. 108.

**More Words.**

At this point I trust that I have presented enough evidence to throw serious doubt on most of the standard stories surrounding the Cardano-Tartaglia affair. The exercise has not been, however, merely to bring to light careless errors in the popular and semi-popular literature. Mistakes, after all, are inevitable. If, however, we extend the concept of scholar to include writers and editors, to any profession that strives toward getting at truth rather than hoodwinking an audience, then it seems to me that such callings require not only intellectual honesty, but intellectual discipline and a basic attention to detail, where the devil resides. The fact that, on the one hand, the Wikipedia editors get the date of Tartaglia's death correct, but on the other hand repeat the story that he abetted Cardano's arrest, tempts one to laugh. As mentioned in the introduction, the apocrypha I've discussed never seem to be repeated in the more scholarly works about Cardano that concern his astrological or medical activities. The tales are apparently confined to the mathematical sphere. Alan Wykes may not have been a mathematician, but many of his readers seem to be. Scientists are suckers.[50]

---

[50] Eckman, p. 77.

Experience forewarns that a non-negligible percentage of readers will meet the present essay with a shrug and reply that legends, at least great ones, are preferable to mundane "true" stories.  A first answer is that, yes, great legends confer moral truths.  In this case of the Great Feud I do not see any deep truths, only negligence, deception and mean-spiritedness.   I can provide a second answer by recounting yet another tale:  Thirty years ago I published an investigation on the various myths surrounding Evariste Galois.[51]  Intending to announce my findings at a seminar at the University of Texas, I thought it would be appropriate to wear a period costume for the occasion and betook myself to the drama department.  The wardrobe mistress didn't have anything on hand from the proper timeframe, and so I asked her just to give me a nice ruffled shirt.   At this she took offense, saying that I was concerned only with historical accuracy in science, not in costumes.  She did relent and lent me a beautiful shirt, but the lesson was a good one and has remained.   Scientists only reluctantly acknowledge truth in other fields, but standards are standards.  If one prefers tall tales and inventions to research, that's fine, but don't call it history.

---

[51] Original version: Tony Rothman, "Genius and Biographers: The Fictionalization of Evariste Galois," Amer. Math. Mon. **89**, 84 (1982).  Revised version available online at various locations.

Acknowledgements

My thanks go to Enrico Lorenzini and Robert Kaster for assistance with translations and to Peter Pesic and the referee for helpful comments.

Tony Rothman's latest books are *Firebird*, a novel about a race for nuclear fusion and *The Course of Fortune*, a historical novel about the Great Siege of Malta.